\documentclass[12pt]{article}
\usepackage[utf8]{inputenc}
\usepackage{amsmath,amssymb,amsthm}
\renewcommand\le{\leqslant}
\renewcommand\ge{\geqslant}

\newtheorem*{theorem}{Theorem}
\newtheorem*{lemma}{Lemma}

\newcommand\eps{\varepsilon}

\title{Kolmogorov widths of the class $W_1^1$}
\author{Yuri Malykhin}

\begin{document}
\maketitle
\begin{abstract}
    We prove that $d_n(W^1_1,L_q)\asymp n^{-1/2}\log n$, $2<q<\infty$. This
    completes the study of orders of decay of Kolmogorov widths for the
    classical case of the univariate
    Sobolev classes of integer smoothness.
\end{abstract}

\paragraph{Introduction.}
Recall that
the Kolmogorov $n$-width of a set $K$ in a normed space $X$ is defined as the
error of the best approximation of $K$ by $n$-dimensional subspaces in $X$:
$$
d_n(K,X) := \inf_{\substack{Q_n\subset X\\\dim Q_n\le n}}\sup_{x\in K}\inf_{y\in Q_n}\|x-y\|_X.
$$
This notion was introduced by Kolmogorov in~\cite{Kol36}.
He considered the Sobolev class $W^r_2$ on the segment $[0,1]$ and expressed
$d_n(W^r_2,L_2)$ in terms of the spectrum of an appropriate differential operator.

The orders of the decay for the widths $d_n(W^r_p[0,1],L_q[0,1])$ were studied
in 1960s and 1970s (Kolmogorov, Stechkin, Tikhomirov, Ismagilov, Kashin, Gluskin,
and others~--- see Chs.13,14 in~\cite{LGM}).
The major breakthrough in this problem was made by Kashin~\cite{Kas77} who
settled the case $q>\max\{2,p\}$;
this paper also contains historical references.
So, sharp order estimates were obtained for almost all of the parameters
$r\in\mathbb N$, $1\le p,q\le\infty$.

Only the following case was left open: $r=p=1$, $2<q<\infty$. The class
$W_1^1:=W^1_1[0,1]$ consists of absolutely
continuous functions $f$ such that $\int_0^1|f'(x)|\,dx\le 1$. Kulanin~\cite{Kul83}
proved the following bounds:
\begin{equation}
    \label{kulanin}
    n^{-1/2}\log^{1/2-\eps}n \lesssim d_n(W^1_1,L_q) \lesssim n^{-1/2}\log n,
\end{equation}
for any fixed $\eps>0$ ($A\lesssim B$ means $A=O(B)$; we write $A\asymp B$ if $A\lesssim B$
and $A\gtrsim B$); hence, the order of the width contains a logarithmic term.
Belinsky~\cite{Bel84} proved that the upper bound is sharp if we consider
only subspaces generated by $n$ arbitrary harmonics:
$d_n^T(W^1_1,L_q) \asymp n^{-1/2}\log n$.

Later, in~\cite{KMR18} the $\eps$ in~\eqref{kulanin} was removed, but the
logarithmic gap remained. The author studied the
closely related case of Besov classes $B^1_{1,\theta}$ and obtained sharp estimates
for their widths~\cite{Mal21}. Unfortunately, this result did not give better
bounds for $W_1^1$.
In this paper we finally remove the logarithmic gap.

\begin{theorem}
    For any $q\in(2,\infty)$ there exist quantities $c_1(q),c_2(q)>0$, such that
    for all $n\ge 2$ the inequality holds
    $$
    c_1(q)n^{-1/2}\log n \le d_n(W_1^1[0,1],L_q[0,1]) \le c_2(q)n^{-1/2}\log n.
    $$
\end{theorem}

In the proof we use standard dyadic decomposition~\eqref{haar},~\eqref{lq_pachki} as
in~\cite{Kul83}. The drawback of his method was that all dyadic
levels were considered independently. Here we control them simultaneously, using
the specially tailored duality and averaging~\eqref{drob1},~\eqref{drob2}, \eqref{isotrop}
in a similar vein as Gluskin~\cite{Gls81}. See also recent
papers~\cite{Mal24,MR25} on the average widths.

Consider the ``step'' functions $\chi_t$ that are defined as $\chi_t(x):=1$,
$0\le x<t$ and $\chi_t(x):=-1$, $t\le x\le 1$.

\begin{lemma}
    Let $q\in(2,\infty)$, $\gamma>1/q$. There exists
    $\delta=\delta(q,\gamma)>0$, such that
    for any sufficiently large $n$
    and any family of functions $\{\eta_t(\cdot)\}_{0\le t\le 1}$ that lie in
    some $n$-dimensional subspace of $L_q[0,1]$, at least one of the following
    inequalities holds:
    \begin{equation}
        \label{l2low}
    \left(\int_0^1 \|\chi_t-\eta_t\|_{L_2}^2\,dt\right)^{1/2} \ge \delta
    n^{-1/2}\log n,
    \end{equation}
    \begin{equation}
        \label{lqlow}
        \left(\int_0^1 \|\chi_t-\eta_t\|_{L_q}^q\,dt\right)^{1/q} \ge n^{-\gamma}.
    \end{equation}
\end{lemma}

Now the proof of Theorem is almost immediate.
\begin{proof}
The upper bound is well-known, see~\cite{Kul83}. Let us prove the lower bound.
Fix any $\gamma\in(1/q,1/2)$ and obtain the corresponding $\delta$ from Lemma.
Suppose that $d_n(W_1^1,L_q) < (\delta/2)n^{-1/2}\log n$.
Functions $\chi_t$ lie in the closure (in $L_q$) of $2W_1^1$;
therefore, there exists a subspace $Q_n$ and functions $\eta_t\in Q_n$, such that
    $\sup_{0\le t\le 1} \|\chi_t-\eta_t\|_{L_q} < \delta n^{-1/2}\log n$.
Then~\eqref{l2low} is clearly violated and hence
$$
\sup_t \|\chi_t - \eta_t\|_{L_q} \ge 
\left(\int_0^1 \|\chi_t - \eta_t\|_{L_q}^q\,dt\right)^{1/q} \ge n^{-\gamma}.
$$
This contradicts our assumption for large $n$ (since $\gamma<1/2$) and establishes the lower bound on
the width.
\end{proof}

\paragraph{Proof of Lemma.}

Suppose that~\eqref{l2low} is violated (we will choose $\delta$ later),
let us prove~\eqref{lqlow}.
Denote by $c_{k,j}(f)$ the coefficients of the expansion of the function
$f$ in the $L_\infty$-normalized Haar system $\{h_{k,j}\}$:
\begin{equation}
    \label{haar}
c_{k,j}(f) := 2^k\int_0^1 f(x)h_{k,j}(x)\,dx
= 2^k(\int_{\frac{j-1}{2^{k}}}^{\frac{j-1/2}{2^{k}}}f(x)\,dx -
\int_{\frac{j-1/2}{2^{k}}}^{\frac{j}{2^{k}}}f(x)\,dx).
\end{equation}
Let $P_kf := \sum_{j=1}^{2^k}c_{k,j}(f)h_{k,j}$ be the $k$-th dyadic pack of the
Fourier--Haar series of $f$.

Fix some numbers $k_0<k_1$ (see~\eqref{k0k1}) and consider the vector space
$$
\mathcal{T} := \{x=(x_{k,j})\colon k_0 \le k\le k_1,\;j=1,\ldots,2^k\}
$$
with weights $w(k,j)=2^{-k}$. In this weighted space the scalar product and the
$\ell_p$-norms are defined:
$$
\langle x,y\rangle_w := \sum_{k=k_0}^{k_1} 2^{-k}\sum_{j=1}^{2^k}x_{k,j}y_{k,j},
\quad
\|x\|_{p,w} := \left(\sum_{k=k_0}^{k_1} 2^{-k}\sum_{j=1}^{2^k}
|x_{k,j}|^p\right)^{1/p}.
$$
We will use that $\langle x,y\rangle_w \le \|x\|_{p,w}\|y\|_{p',w}$.

Well-known Marcinkiewicz--Paley theorem gives (note that $q>2$):
$$
\|f\|_{L_q}^q\asymp \|((\int_0^1 f(x)\,dx)^2 + \sum_{k\ge 0} P_k^2f)^{1/2}\|_{L_q}^q
\ge \sum_{k\ge 0} \|P_kf\|_{L_q}^q
\ge \sum_{k=k_0}^{k_1} \|P_kf\|_{L_q}^q,
$$
\begin{equation}
    \label{lq_pachki}
    \|f\|_{L_q}^q \gtrsim \sum_{k=k_0}^{k_1}2^{-k}\sum_{j=1}^{2^k}|c_{k,j}(f)|^q
    = \|c_{k,j}(f)\|_{q,w}^q.
\end{equation}
We also use the simple inequality
\begin{equation}
    \label{l2_pachki}
    \|f\|_{L_2}^2 \ge \sum_{k=k_0}^{k_1} \|P_kf\|_{L_2}^2 =
    \|c_{k,j}(f)\|_{2,w}^2.
\end{equation}

Define
$$
X_{k,j}(t) := c_{k,j}(\chi_t),\quad H_{k,j}(t) := c_{k,j}(\eta_t).
$$
We shall consider $X$, $H$ as the collection of functions of a variable
$t\in[0,1]$ and at the same time as the parametrized family of vectors
$X(t),H(t)\in \mathcal{T}$.
It is rather obvious that $X_{0,1}(t)=1-|2t-1|$; other functions are also
piecewise--linear and may be expressed as
$X_{k,j}(t) = X_{0,1}(2^kt-(j-1))$ for $t\in[(j-1)2^{-k},j2^{-k}]$ and
$X_{k,j}(t)=0$ outside this interval.

Vectors $H(t)$ lie in some $n$-dimensional subspace
$V_n\subset \mathcal{T}$ for all $t$, because $\eta_t$ lie in a $n$-dimensional
subspace and the map $f\mapsto (c_{k,j}(f))$ is linear.

Let $Z_{k,j}(t)$ be some functions. One can use the duality and the Holder
inequality and obtain:
\begin{multline}\label{drob1}
    \int_0^1 \langle X(t)-H(t),Z(t)\rangle_w\,dt
    \le \int_0^1 \|X(t)-H(t)\|_{q,w}\cdot \|Z(t)\|_{q',w}\,dt \le \\
    \le (\int_0^1\|X(t)-H(t)\|_{q,w}^q\,dt)^{1/q}(\int_0^1\|Z(t)\|_{q',w}^{q'}\,dt)^{1/q'}.
\end{multline}
Therefore, using~\eqref{lq_pachki}, we have
\begin{equation}\label{drob2}
    (\int_0^1 \|\chi_t-\eta_t\|_q^q\,dt)^{1/q} \gtrsim (\int_0^1\|X(t)-H(t)\|_{q,w}^q\,dt)^{1/q} \ge \frac{I_1 - I_2}{I_3^{1/q'}},
\end{equation}
where $I_1 := \int_0^1\langle X(t),Z(t)\rangle_w\,dt$,
$I_2 := \int_0^1\langle H(t),Z(t)\rangle_w\,dt$,
$I_3 := \int_0^1\|Z(t)\|_{q',w}^{q'}\,dt$.

Now let us make our construction specific in order to estimate
$I_1,I_2,I_3$:
\begin{equation}
    \label{k0k1}
    k_0 := \lfloor \log_2 n\rfloor,
    \quad k_1 := \lfloor \gamma q\log_2 n\rfloor,
\end{equation}
$$
Z_{k,j}(t) := \begin{cases}
    2^{k}a(X_{k,j}(t)-1/2),&\quad x\in[(j-1)2^{-k},j2^{-k}],\\
    0,&\quad \mbox{otherwise},
\end{cases}
$$
where $a:=\sqrt{12}$.

Note that $\int_0^1Z_{k,j}(t)\,dt=0$ and
$\int_0^1(t-t_0)Z_{k,j}(t)\,dt=0$ for any $t_0$ (it is enough to consider
$t_0$ in the middle of the support of $Z_{k,j}$).
Let us prove that
\begin{equation}
    \label{orthonormal}
\int_0^1 Z_{k,j}(t)Z_{k',j'}(t)\,dt=0,
\quad \int_0^1 Z_{k,j}(t)^2\,dt = 2^{k},
    \end{equation}
for all $(k,j)\ne (k',j')$. First we check the orthogonality. If $k=k'$, the
supports of functions do not intersect.
If $k<k'$ and the supports intersect, then on the support of the function
$Z_{k',j'}$ the function $Z_{k,j}$ is linear and the integral of their product
is zero. Moreover,
$$
\int_0^1 Z_{k,j}^2(t)\,dt = 2^{2k}a^2
\int_{\frac{j-1}{2^k}}^{\frac{j}{2^k}}
(X_{k,j}(t)-1/2)^2\,dt = 2^k a^2\int_0^{1}(X_{0,1}(t)-1/2)^2\,dt = 2^k.
$$

Proceed to the quantities $I_1,I_2,I_3$. First we compute $I_1$:
$$
\int_0^1 X_{k,j}(t)Z_{k,j}(t)\,dt
= 2^k a\int_{\frac{j-1}{2^{k}}}^{\frac{j}{2^k}}(X_{k,j}(t)-1/2)^2\,dt =
a\int_0^1(X_{0,1}(t)-1/2)^2\,dt = 1/a,
$$
$$
I_1 = \int_0^1\langle X(t),Z(t)\rangle_w\,dt
= \sum_{k=k_0}^{k_1} 2^{-k} \sum_{j=1}^{2^k} \int_0^1
X_{k,j}(t)Z_{k,j}(t)\,dt
= (k_1-k_0+1)/a.
$$

Equalities~\eqref{orthonormal} imply that $Z$ is isotropic in $\mathcal T$, that
is $\int_0^1 \langle v,Z(t)\rangle_w^2\,dt = \langle v,v\rangle_w^2$ for any
$v\in\mathcal T$. 
Let $\Pi_{V_n}$ be the orthoprojector in $\mathcal T$ on the subspace $V_n$
(that contains vectors $H(t)$). Due to the isotropy we have (where $v_1,\ldots,v_n$
is an orthonormal basis of $V_n$):
\begin{equation}
    \label{isotrop}
\int_0^1\|\Pi_{V_n}Z(t)\|_{2,w}^2\,dt = \int_0^1 \sum_{m=1}^n\langle
v_m,Z(t)\rangle_w^2\,dt = n.
\end{equation}
This allows us to estimate $I_2$:
\begin{multline*}
I_2 = \int_0^1\langle H(t),Z(t)\rangle_w\,dt
= \int_0^1\langle H(t),\Pi_{V_n}Z(t)\rangle_w\,dt \le \\
\le (\int_0^1\|H(t)\|_{2,w}^2\,dt)^{1/2}(\int_0^1 \|\Pi_{V_n}Z(t)\|_{2,w}^2\,dt)^{1/2}
    = n^{1/2} (\int_0^1\|H(t)\|_{2,w}^2\,dt)^{1/2}.
\end{multline*}
The inequality~\eqref{l2_pachki} gives:
\begin{multline*}
    \|H(t)\|_{2,w}^2 \le 2\|X(t)\|_{2,w}^2 + 2\|X(t)-H(t)\|_{2,w}^2
    \le 2\|X(t)\|_{2,w}^2 + 2\|\chi_t-\eta_t\|_{L_2}^2.
\end{multline*}
The first term is bounded for all $t$:
$\|X(t)\|_{2,w}^2\le \sum_{k\ge k_0} 2^{-k} \le 2^{1-k_0}\le 4/n$.
The second term is bounded in average, since~\eqref{l2low} is violated, so
$$
\int_0^1\|H(t)\|_{2,w}^2\,dt \lesssim n^{-1} + 
\int_0^1 \|\chi_t-\eta_t\|_{L_2}^2\,dt \lesssim \delta^{2} n^{-1}\log^2 n.
$$
Hence, $I_2 \lesssim \delta \log n$ for large $n$.
Finally,
\begin{multline*}
I_3 = \int_0^1\|Z(t)\|_{q',w}^{q'}\,dt
= \sum_{k=k_0}^{k_1} 2^{-k} \sum_{j=1}^{2^k} 2^{kq'}a^{q'}
    \int_{\frac{j-1}{2^k}}^{\frac{j}{2^k}}|X_{k,j}(t)-1/2|^{q'}\,dt \lesssim \\
\lesssim \sum_{k=k_0}^{k_1} 2^{k(q'-1)} \asymp 2^{k_1(q'-1)}.
\end{multline*}
So, $I_3^{1/q'} \lesssim 2^{k_1(q'-1)/q'} = 2^{k_1/q} \le n^\gamma$.

The numerator in~\eqref{drob2} is at least $(k_1-k_0+1)/a - C \delta \log n \ge
((\gamma q-1)/a-C\delta)\log n \gtrsim \log n$, is  $\delta$ is small enough.
The denonimator is $\lesssim n^{\gamma}$. Q.E.D.

\paragraph{Acknowledgements.} 
Author thanks K.S.~Ryutin for fruitful discussions.

\end{document}